\documentclass{amsart}
\usepackage{amsmath}
\usepackage{amscd}
\usepackage{amssymb}

\newtheorem{theorem}{Theorem}
\newtheorem{corollary}[theorem]{Corollary}

\newtheorem{proposition}{Proposition}
\newtheorem{definition}{Definition}

\newtheorem*{ackn}{Acknowledgements}

\newcommand{\supp}{\operatorname{supp}}

\newcommand{\sign}{\operatorname{sign}}

\newcommand{\N}{\mathbb{N}}

\newcommand{\Hb}{\mathcal{H}}

\title[Embeddings into $\ell_p$-spaces]
{\large On coarse embeddability into $\ell_p$-spaces and a
conjecture of Dranishnikov}

\author{Piotr W. Nowak}

\address{Department of Mathematics,
Vanderbilt University, 1326 Stevenson Center, Nashville, TN 37240
USA.}

\subjclass[2000]{Primary 46C05; Secondary 46T99}

\email{pnowak@math.vanderbilt.edu}

\keywords{coarse embeddings, $\ell_p$-spaces, property A, Novikov
Conjecture}

\begin{document}

\begin{abstract}
 We show that the Hilbert space is coarsely embeddable into any $\ell_p$
for $1\le p<\infty$. In particular, this yields new
characterizations of embeddability of separable metric spaces into
the Hilbert space.
\end{abstract}
 \maketitle
 Coarse embeddings were
defined by M.~Gromov \cite{gromov} to express the idea of inclusion
in the large scale geometry of groups. G.~Yu showed later that the
case when a finitely generated group with a word length metric is
being embedded into the Hilbert space is of great importance in
solving the Novikov Conjecture \cite{yu-embeddings}, while recent
work of G.~Kasparov and G.~Yu \cite{kasparov-yu} treats the case
when the Hilbert space is replaced with just a uniformly convex
Banach space. Due to these remarkable theorems coarse embeddings
gain a great deal of attention, but still embeddability into the
Hilbert, and more generally Banach spaces, is not entirely
understood with many question remaining open.

In this context the class of $\ell_p$-spaces seems to be
particularly interesting. Their embeddability into the Hilbert space
is known - $\ell_p$ admits such an embedding when $0<p\le 2$ but do
not if $p>2$ due to a recent result of W.~Johnson and
N.~Randrianarivony \cite{johnson-randrianarivony}. In this note we
study the opposite situation, i.e. we show that the separable
Hilbert space embeds into $\ell_p$ for any $1\le p <\infty$. As a
consequence we obtain a new characterization of embeddability into
$\ell_2$, namely that embeddings into $\ell_p$ for $1\le p \le 2$
are all equivalent.

In \cite[Section 6]{guentner-kaminker} the authors advertised a
conjecture  stated by A.N.~Dranishnikov \cite[Conjecture
4.4]{dranishnikov}: \emph{a discrete metric space has property A if
and only if it admits a coarse embedding into the the space
$\ell_1$}. The results presented in this note show, that this is the
same as asking whether property A is equivalent to embeddability
into the Hilbert space, and although it is a folk conjecture that
such statement is not true, no example distinguishing between the
two is known.

\begin{ackn}\normalfont
I would like to thank Guoliang Yu for inspiring conversations on
coarse geometry of Banach spaces.\\
\end{ackn}
\section*{$L_p$-spaces and the Mazur Map}
In everything what follows we consider only separable
$L_p(\mu)$-spaces and we use the standard notation
$\ell_p=\ell_p(\N)$. By $S(X)$ we denote the unit sphere in the
Banach space $X$.

The Mazur map $M_{p,q}:S(\ell_p)\to S(\ell_q)$ is defined by the
formula
$$M_{p,q}(x)=\left\{\vert x_i\vert^{\frac{p}{q}} \sign{x_i}\right\}_{i=1}^{\infty}$$
where $x=\lbrace x_i\rbrace_{i=1}^{\infty}\in \ell_p$. It is a
uniform homeomorphism between unit spheres of $\ell_p$-spaces. More
precisely, it satisfies the following inequalities:

\begin{equation}\label{mazur.map.estimates}
\frac{p}{q}\Vert x-y\Vert_p\ \le\ \Vert
M_{p,q}(x)-M_{p,q}(y)\Vert_q\ \le\ C\,\Vert x-y\Vert_p^{p/q}
\end{equation}
for all $x,y\in S(\ell_p)$ and $p<q$ and the opposite inequalities
if $p>q$ (note that $M_{p,q}=M_{q,p}^{-1}$), where the constant $C$
depends only on $\frac{p}{q}$. For the proof of these estimates and
details on the Mazur map and its applications we refer the reader to
\cite[Chapter 9.1]{bl}.

If $\left\{X_n\right\}_{n\in\N}$ is a sequence of Banach spaces,
we denote by $\left(\sum X_n \right)_p$ the direct sum of $X_n$
with the $p$-norm, i.e.$$\left(\sum_{n=1}^{\infty} X_n \right)_p =
\left\{ \mathbf{x}=\{x_n\}_{n\in\N} \ :\ x_n\in X_n,
\sum_{n=1}^{\infty}\Vert x_n\Vert^p <\infty \right\},$$
$$\Vert\mathbf{x}\Vert_p=\left(\sum_{n=1}^{\infty}\Vert
x_n\Vert^p\right)^{\frac{1}{p}}.$$
 Recall that $\ell_p$ is isometric to $\left(\sum
\ell_p
\right)_p$\\

We will also need the following classification of separable
$L_p$-spaces.
\begin{theorem}[see e.g. {\cite[III.A]{wojtaszczyk}}]\label{class.Lp}
A separable space $L_p(\mu)$ is isometric to one of the following
spaces: $\ell_p^n$ for $n=1,2,3,...$, $L_p[0,1]$, $\ell_p$,
$(L_p[0,1]\oplus\ell_p^n)_p$ for $\ n=1,2,3,...$,
$(L_p[0,1]\oplus\ell_p)_p$.\\
\end{theorem}

\section*{A condition for coarse embeddability}
We recall the definition of a coarse embedding.
\begin{definition}\label{defnition - coarse embedding}
\normalfont Let $X,Y$ be metric spaces. A map $f\colon X\to Y$ is
a \emph{coarse embedding} if there exist non-decreasing functions
$\rho_1,\rho_2\colon [0,\infty)\to[0,\infty)$ satisfying
\begin{enumerate}
\item $\rho_1(d_X(x,y))\le d_Y(f(x),f(y))\le \rho_2(d_X(x,y))$\; for all \;$x,y\in X$,
\item $\lim_{t\to\infty}\rho_1(t)=+\infty$.
\end{enumerate}
\end{definition}

In \cite{dg} M.~Dadarlat and E.~Guentner characterized spaces
coarsely embeddable into the Hilbert $\Hb$ space in terms of
existence of maps into the unit sphere $S(\Hb)$. Their result is a
reminiscence of a characterization of \emph{uniform embeddability}
(meaning an existence of a uniform homeomorphism onto a subset)
into a Hilbert space obtained by Aharoni et al in \cite{amm}.

\begin{theorem}[{\cite[Theorem 2.1]{dg}}]\label{dg.embeddability.condition}
A metric space $X$ admits a coarse embedding into the Hilbert
space $\Hb$ if and only if for every $R>0$ and $\varepsilon>0$
there is a map $\varphi :X\to S(\Hb)$ and $S>0$ satisfying
\begin{enumerate}
\item $\sup \left\{ \Vert \varphi(x)-\varphi(y)\Vert_{\Hb}\ :\ x,y\in
X,\ d(x,y)\le R\right\}\le \varepsilon,$
\item $\lim_{S\to\infty} \inf \left\{ \Vert \varphi(x)-\varphi(y)\Vert_{\Hb}\
:\ x,y\in X,\ d(x,y)\ge S\right\}=2$
\end{enumerate}
\end{theorem}

 We are going to use this idea
to prove a similar condition for embeddings into the spaces
$\ell_p$. The proof relies on the original proof of Theorem
\ref{dg.embeddability.condition}.
\begin{theorem}\label{weakA}
Let $X$ be a metric space and $1\le p <\infty$. If there is a
$\delta>0$ such that for every $R>0$, $\varepsilon>0$ there is a map
$\varphi :X\to S(\ell_p)$ satisfying
\begin{enumerate}
\item $\sup \left\{ \Vert \varphi(x)-\varphi(y)\Vert_p\ :\ x,y\in
X,\ d(x,y)\le R\right\}\le \varepsilon,$
\item $\lim_{S\to\infty} \inf \left\{ \Vert \varphi(x)-\varphi(y)\Vert_p\
:\ x,y\in X,\ d(x,y)\ge S\right\}\ge\delta,$
\end{enumerate}
then $X$ admits a coarse embedding into $\ell_p$.
\end{theorem}

\begin{proof}
By the assumptions for every $n\in\N$ there is a map
$\varphi_n:X\to S(\ell_p)$ and a number $S_n>0$ such that
$\Vert\varphi_n(x)-\varphi_n(y)\Vert_p\le\frac{1}{2^n}$ whenever
$d(x,y)\le n$ and $\Vert\varphi_n(x)-\varphi_n(y)\Vert_p\ge
\frac{\delta}{2}$ whenever $d(x,y)\ge S_n$. Without loss of
generality we can choose the sequence of $S_n$'s to be strictly
increasing and tending to infinity as $n\to\infty$.

Choose $x_0\in X$ and define a map $\Phi:X\to \left(\sum
\ell_p\right)_p$ by the formula
$$\Phi(x)= \oplus_{n=1}^{\infty} \left(\varphi_n(x)-\varphi_n(x_0)\right).$$
It is easy to see that
$$\Vert \Phi(x)\Vert_p^p= \sum_{n=1}^{\infty}
\left\Vert \varphi_n(x)-\varphi_n(x_0)\right\Vert_p^p<\infty$$
which shows that $\Phi$ is well-defined.

We will show that $\Phi$ is a coarse embedding. Take $k\in \N$ and
$\sqrt[p]{k-1}\le d(x,y)<\sqrt[p]{k}$. Then
$$\Vert \Phi(x)-\Phi(y)\Vert_p^p=\sum_{n=1}^{k-1}\left\Vert
\varphi_n(x)-\varphi_n(y)\right\Vert_p^p+\sum_{n=k}^{\infty}\left\Vert
\varphi_n(x)-\varphi_n(y)\right\Vert_p^p$$ $$\le
2^p(k-1)+\sum_{n=k}^{\infty} \frac{1}{2^{kp}}\le 2^p(k-1)+1\le 2^p
d(x,y)^p+1$$ The first estimate comes from the fact that unit
vectors cannot be more than distance 2 apart.

On the other hand for $S_{k-1}\le d(x,y)< S_k$ we have
$$\Vert \Phi(x)-\Phi(y)\Vert_p^p\ge\sum_{n=1}^{k-1}\left\Vert \varphi_n(x)-\varphi_n(y)\right\Vert_p^p\ge (k-1)\left(\frac{\delta}{2}\right)^p.$$
Thus we can choose
$\rho_1(t)=\sum_{n=1}^{\infty}\,\,\delta\sqrt[p]{n}
\chi_{[S_{n-1},S_n)}(t)$, $\rho_2(t)=2t+1$ and it is clear that
$\Phi$ is a coarse embedding.
\end{proof}

G.~Yu defined property A \cite{yu-embeddings}, which gives a
sufficient condition for embeddability of a discrete metric space
into a Hilbert space. We recall a characterization of property A
given by J.L.~Tu.

\begin{proposition}[\cite{tu}]
A metric space $X$ has property A if and only if for every $R>0$ and
$\varepsilon>0$ there is a map $\eta:X\to S(\ell_2(X))$ and $S>0$
such that
\begin{enumerate}
\item $\Vert \eta(x)-\eta(y)\Vert_2\le\varepsilon$ when $d(x,y)\le
R$;
\item $\supp \eta(x) \subset B(x,S)$ for all $x\in X$.
\end{enumerate}
\end{proposition}
Theorem \ref{dg.embeddability.condition} and the above
characterization exhibit the subtle relation between Property A and
coarse embeddability.\\

The following proposition shows that the property of Theorem
\ref{weakA} is not sensitive to changing the index $p$.

\begin{proposition}\label{weakA.for.p.then.q}
Let $X$ have the property described in Theorem \ref{weakA} with
respect to some $1\le p <\infty$. Then $X$ has the same property
with respect to any $1\le q<\infty$
\end{proposition}
\begin{proof}
For $R>0$ and $\varepsilon>0$ given a map $f_p:X\to S(\ell_p)$ which
satisfies conditions (1) and (2) of Theorem \ref{weakA} define
$f_q:X\to S(\ell_q)$ by the formula
$$ f_q(x)=M_{p,q}\left[f_p(x)\right],$$
where $M_{p,q}:S(\ell_p)\to S(\ell_q)$ is the Mazur map.

If $p<q$, by inequalities (\ref{mazur.map.estimates}) we have
$$\frac{p}{q}\Vert f_p(x)-f_p(y)\Vert_p \le \Vert
f_q(x)-f_q(y)\Vert_q\le C \Vert f_p(x)-f_p(y)\Vert_p^{p/q}.$$
Consequently
$$\sup \left\{ \Vert f_q(x)-f_q(y)\Vert_q\ :\ x,y\in X,\
d(x,y)\le R\right\}\le C \varepsilon^{p/q},$$ and
$$\lim_{S\to\infty} \inf
\left\{ \Vert f_q(x)-f_q(y)\Vert_q\ :\ x,y\in X,\ d(x,y)\ge
S\right\}\ge\frac{p}{q}\delta.$$ The case $p>q$ is proved
similarly.
\end{proof}

In the case of Property A a statement similar to Proposition
\ref{weakA.for.p.then.q} was studied by Dranishnikov under the name
of Property A$_p$ in \cite{dranishnikov}.

\begin{corollary}
If $X$ admits a coarse embedding into $\ell_2$ then it admits a
coarse embedding into any $\ell_p$ with $1\le p <\infty$. In
particular, the separable Hilbert space embeds into all $\ell_p$.
\end{corollary}

\begin{proof}
If $X$ admits a coarse embedding into $\ell_2$ then, by Theorem
\ref{dg.embeddability.condition}, $X$ has the property from
Proposition \ref{dg.embeddability.condition} for $\ell_2$. By
Proposition \ref{weakA.for.p.then.q} has this property also for
$\ell_p$, $1\le p<\infty$ and by Proposition \ref{weakA} admits an
embedding into $\ell_p$.
\end{proof}

In \cite{nowak} it was shown that $L_p(\mu)$ for $1\le p\le 2$ admit
a coarse embedding into the Hilbert space and that coarse
embeddability into $\ell_2$ is equivalent to coarse embeddability
into $L_p[0,1]$ again for $1\le p\le 2$. This allows us to state

\begin{theorem}
Let $X$ be a separable metric space. Then the following conditions
are equivalent:
\begin{enumerate}
\item $X$ admits a coarse embedding into the Hilbert space;
\item $X$ admits a coarse embedding into $\ell_p$ for any $1\le p
<2$;
\item $X$ admits a coarse embedding into $L_p[0,1]$ for any $1\le
p<2$.
\end{enumerate}
\end{theorem}

Note that this covers all separable $L_p(\mu)$-spaces with $1\le
p\le 2$. This is particularly interesting since the spaces $L_p$ for
different $p$'s are not coarsely equivalent. To see this assume they
are and take $f:L_p(\mu)\to L_q(\mu)$ to be the coarse equivalence.
Since $L_p$-spaces are geodesic, $f$ is in fact a quasi-isometry and
it induces a Lipschitz equivalence on their ultrapowers. By a
theorem of Heinrich \cite{heinrich} ultrapowers of $L_p$ spaces are
again $L_p$ spaces (possibly on a different measure), and the
assertion follows from a classical fact that Lipschitz equivalence
on $L_p$-spaces induces a linear isomorphism.

\bibliographystyle{amsalpha}

\end{document}